\documentclass[12pt]{amsart}
\usepackage{amssymb,latexsym,eufrak,amsmath,amscd, graphicx}

\usepackage{amsfonts}
\usepackage{amscd}

\def\ZZ{{\mathbb Z}}
\def\NN{{\mathbb N}}

\def\CC{{\mathbb C}}
\def\AA{{\mathbb A}}
\def\RR{{\mathbb R}}
\def\QQ{{\mathbb Q}}

\DeclareMathOperator{\ord}{ord}

\newtheorem{lemma}{Lemma}[section]
\newtheorem{theorem}[lemma]{Theorem}

\newtheorem{proposition}[lemma]{Proposition}

\newtheorem{conjecture}[lemma]{Conjecture}

\theoremstyle{definition}

\newtheorem{remark}[lemma]{Remark}

\theoremstyle{remark}
\newtheorem*{remark*}{Remark}
\newtheorem*{note*}{Note}

\begin{document}

\title{On Igusa zeta functions of monomial ideals}

\author[J. Howald]{Jason Howald}
\address{ (Howald) Department of Mathematics and Computer Science,
John Carroll University, 20700 North Park Blvd., University Heights,
OH 44118, USA} \curraddr{SUNY Potsdam, Dept. of Mathematics, 44
Pierrepont Avenue, Potsdam, NY 13676-2294, USA}
 \email{{\tt howaldja@potsdam.edu }}

\author[M. Musta\c{t}\v{a}]{Mircea~Musta\c{t}\v{a}}
\address{(Musta\c{t}\v{a}) Department of Mathematics, University of Michigan,
Ann Arbor, MI 48109, USA}
\email{{\tt mmustata@umich.edu}}

\author[C. Yuen]{Cornelia Yuen}
\address{(Yuen) Department of Mathematics, University of Michigan,
Ann Arbor, MI 48109, USA} \curraddr{Department of Mathematics,
University of Kentucky, 825 Patterson Office Tower, Lexington, KY
40506, USA}

\email{{\tt cyuen@ms.uky.edu  }}

\begin{abstract}
We show that the real parts of the poles of the Igusa zeta function
of a monomial ideal can be computed from  the torus-invariant
divisors on the normalized blow-up of the affine space along the
ideal. Moreover, we show that every such number is a root of the
Bernstein-Sato polynomial associated to the monomial ideal.
\end{abstract}

\thanks{2000\,\emph{Mathematics Subject Classification}.
Primary 14B05; Secondary 14M25}
\thanks{ The
research of the second author was partially supported by NSF grant
DMS 0500127.}

\maketitle

\markboth{J. HOWALD, M. MUSTA\c T\v A and C.~ YUEN} {IGUSA ZETA
FUNCTIONS OF MONOMIAL IDEALS}

\section{Introduction}

If $f$ is a nonconstant polynomial in $\ZZ[x_1,\ldots,x_n]$ and $p$
is a fixed prime, then the Igusa zeta function of $f$ is defined by
\begin{equation}\label{def1}
Z_f(s)=\int_{(\ZZ_p)^n}|f(y)|_p^s\,|dy|,
\end{equation}
for every $s\in\CC$ with ${\rm Re}(s)>0$. Here $\ZZ_p$ denotes the
ring of $p$-adic integers, endowed with the discrete valuation ${\rm
ord}_p$ and with the $p$-adic absolute value defined by
$$|a|_p=\left(\frac{1}{p}\right)^{{\rm ord}_p(a)}.$$
The measure $\mu$ on $(\ZZ_p)^n$ that is used  in the above integral
is the Haar measure characterized by
$$\mu\left(\prod_{i=1}^np^{m_i}\ZZ_p\right)
=\left(\frac{1}{p}\right)^{\sum_im_i}.$$

As defined, $Z_f$ is a holomorphic function on the half plane
$\{s\mid {\rm Re}(s)>0\}$ and one can show that it admits a
meromorphic extension to $\CC$. In fact, $Z_f$ is a rational
function of $\left(\frac{1}{p}\right)^s$. A proof of rationality
that gives also information on the real parts of the possible poles
of $Z_f$ proceeds as follows. Let $\pi\colon Y\to X=\AA^n$ be an
embedded resolution of singularities of $f$ defined over $\QQ$. This
means that $\pi$ is proper and birational, $Y$ is smooth and the
union of $\pi^*({\rm div}(f))$ and of the exceptional locus of $\pi$
is a divisor with simple normal crossings. For every prime divisor
$E$ in this union, denote by $a_E(f)$ the order of $E$ in
$\pi^*({\rm div}(f))$ and by $k_E$ the order of $E$ in the relative
canonical class $K_{Y/X}$ (this is the divisor locally defined by
${\rm det}({\rm Jac}(\pi))$). Using the Change of Variable formula
for $p$-adic integrals to express $Z_f$ as an integral over
$Y(\ZZ_p)$, Igusa obtained the rationality of $Z_f$ as a function of
$\left(\frac{1}{p}\right)^s$ and the fact that if $s$ is a pole of
$Z_f$, then ${\rm Re}(s)=-\frac{k_E+1}{a_E(f)}$ for some divisor $E$
as above. Our main reference for Igusa zeta functions is Igusa's
book \cite{Ig} (see also Denef's Bourbaki report \cite{De}).

While every divisor on a log resolution of $f$ gives a candidate for
the real part of a pole of $Z_f$, examples show that most of these
numbers do not come actually from poles of $Z_f$. In fact an outstanding open
problem in the field is the following conjecture of Igusa, relating
the poles to one of the basic invariants of the singularities of
$f$, its Bernstein-Sato polynomial.

\begin{conjecture}\label{conj1}
Let $f$ be a non-constant polynomial in $\ZZ[X_1,\ldots,X_n]$. For
almost all primes $p$ the following holds: if $s$ is a pole of
$Z_f$, then the real part of $s$ is a root of the Bernstein-Sato polynomial
$b_f$ of $f$.
\end{conjecture}

We recall that the Bernstein-Sato polynomial of $f$ is a polynomial
in one variable introduced independently in \cite{Be} and \cite{SS}.
It is a subtle but very fundamental invariant of the singularities
of $f$. We do not give its definition as we will not need it, but we
mention that its roots are related to the eigenvalues of the
monodromy of the hypersurface $f^{-1}(0)$. There is, in fact, a
weaker version of the above conjecture that is stated in terms of
these eigenvalues and that is known as the Monodromy Conjecture (see
\cite{De} for more on these conjectures and also \cite{Ve} for some
recent work in this direction).

Our goal in  this note is to prove the analogue of
Conjecture~\ref{conj1} when we replace $f$ by a monomial ideal.
Though less studied, Igusa zeta functions for non-necessarily
principal ideals in $\ZZ[x_1,\ldots,x_n]$ can be defined in a very
similar way with (\ref{def1}). More precisely, if $I$ is a nonzero
proper ideal of $\ZZ[x_1,\ldots,x_n]$ and if we put for $y\in
(\ZZ_p)^n$
 $$\ord_pI(y)=\min\{{\rm ord}_p(f(y))\vert f\in I\},$$
then we have
\begin{equation}
Z_I(s)=\int_{(\ZZ_p)^n}\left(\frac{1}{p}\right)^{\ord_pI(y)}|dy|.
\end{equation}

The above-mentioned results in the case of one polynomial extend in
a straightforward way to the case of an arbitrary ideal. Note that
in order to prove rationality, we need to consider a log resolution
of $I$: this is a morphism $\pi\colon Y\to {\mathbb A}^n$ as before,
such that $\pi^{-1}(V(I))$ is a Cartier divisor and its union with
the exceptional locus of $\pi$ is a divisor with simple normal
crossings. If $E$ is a prime divisor on $Y$ contained in this union,
then $a_E(I)$ is by
definition the coefficient of $E$ in $\pi^{-1}(V(I))$. As in the
case of a principal ideal, one can show that given a log resolution
$\pi$, for every pole $s$ of $Z_I$ there is a divisor $E$ on $Y$
such that ${\rm Re}(s)=-\frac{k_E+1}{a_E(I)}$.

On the other hand, the
definition of the Bernstein-Sato polynomial has been extended in
\cite{BMS3} from the case of one polynomial to that of an arbitrary
ideal. This is again a polynomial in one variable and therefore the
analogue of Conjecture~\ref{conj1} makes sense in this case. We will
prove the monomial case, i.e. when $I$ is generated by monomials.

\begin{theorem}\label{thm1}
If $I$ is a nonzero proper monomial ideal of $\ZZ[x_1,\ldots,x_n]$,
then for every prime $p$ and every pole $s$ of $Z_I$, the real part
of $s$ is a root of the Bernstein-Sato polynomial of $I$.
\end{theorem}

The key ingredient in the proof of the above theorem is a result on
the poles of Igusa-type zeta functions associated to cones. Suppose
that $N\simeq\ZZ^n$ is a lattice and that $\sigma$ is a pointed,
rational, polyhedral cone in $N_{\RR}=N\otimes_{\ZZ}{\RR}$. We denote
by $\sigma^{\circ}$ the relative interior of the cone $\sigma$. If
$\sigma^{\vee}$ is the dual cone in $M_{\RR}$, where $M={\rm
Hom}_{\ZZ}(N,\ZZ)$, and if $\ell_1,\ell_2$ are elements in
$\sigma^{\vee}\cap M$ such that
$\sigma\cap\{v\mid\ell_2(v)=0\}=\{0\}$, then we put
\begin{equation}
Z_{\sigma,\ell_1,\ell_2}(s):=\sum_{v\in\sigma^{\circ}\cap
N}\left(\frac{1}{p}\right)^{\ell_1(v)s+\ell_2(v)}.
\end{equation}
It is easy to see (and it will follow from our computations) that
this is well-defined  and holomorphic in $\{s\in\CC\mid {\rm
Re}(s)>0\}$. We prove the following

\begin{theorem}\label{key}
For every $\sigma$, $\ell_1$ and $\ell_2$ as above,
$Z_{\sigma,\ell_1,\ell_2}$ is a rational function of
$\left(\frac{1}{p}\right)^s$, and therefore can be meromorphically
extended to $\CC$. Moreover, for every pole $s$ of
$Z_{\sigma,\ell_1,\ell_2}$ there is a primitive generator $v$ of a
ray of $\sigma$ such that ${\rm
Re}(s)=-\frac{\ell_2(v)}{\ell_1(v)}$.
\end{theorem}

Given a monomial ideal $I$, we give in the next section a formula
for $Z_I$ in terms of suitable zeta functions for the cones in the
normal fan to the Newton polyhedron of $I$ (we refer for the
relevant definition and for the precise formula to that section).
Let us just mention that this fan defines the toric variety that is
the normalized blowing-up of $\AA^n$ along the ideal $I$. Using this
formula and Theorem~\ref{key}, we will show in \S 3 that the real
part of every pole of $Z_I$ corresponds to a torus-invariant divisor in
the normalized blow-up of $\AA^n$ along $I$ (despite the fact that
the normalized blowing-up is \emph{not} a log resolution of I).

\begin{theorem}\label{thm2}
Let $I$ be a nonzero proper monomial ideal of $\ZZ[x_1,\ldots,x_n]$.
For every pole $s$ of $Z_I$, there is a torus-invariant divisor $E$ on
the normalized blowing-up of $\AA^n$ along $I$ such that
$${\rm Re}(s)=-\frac{k_E+1}{a_E(I)}.$$
\end{theorem}

On the other hand, explicit descriptions of the roots of the
Bernstein-Sato polynomial of a monomial ideal have been obtained in
\cite{BMS1} and \cite{BMS2}. We use the description in \cite{BMS2} and
Theorem~\ref{thm2} to prove Theorem~\ref{thm1} in the last section.

\smallskip

We mention that a description for the Igusa zeta function of a
monomial ideal has also been obtained by Z\'{u}${\rm
\tilde{n}}$iga-Galindo in \cite{Zu}. Moreover, similar results
appear in the work of Denef and Hoornert \cite{DH}, in which one
describes the poles of the Igusa zeta functions for nondegenerate
hypersurfaces with respect to their Newton polyhedron. It is shown
in \emph{loc. cit.} that for such $f$ the real part of essentially
any pole corresponds to a facet of the Newton polyhedron of $f$, as
above. Moreover, Loeser \cite{Lo} showed that under some mild extra
assumptions, these numbers are roots of the Bernstein-Sato
polynomial of $f$, thus proving Conjecture~\ref{conj1} for such
nondegenerate hypersurfaces. On the other hand, note that the
relations between the respective Igusa zeta functions and
Bernstein-Sato polynomials of $f$ and of the corresponding monomial
ideal are not clear in general.

\section{Igusa zeta functions of monomial ideals}

Let $I$ be a nonzero proper ideal of $\ZZ[x_1,\ldots,x_n]$ generated
by monomials. If $u=(u_1,\ldots,u_n)\in\NN^n$, we denote by $x^u$
the corresponding monomial $x_1^{u_1}\ldots x_n^{u_n}$. The Newton
polyhedron $P_I$ of $I$ is the convex hull of those $u$ in $\NN^n$
such that $x^u$ is in $I$.

If $N=\ZZ^n$ and $N_{\RR}=N\otimes_{\ZZ}{\RR}$, we think of $P_I$
as lying in $M_{\RR}$, where $M={\rm Hom}_{\ZZ}(N,\ZZ)$. We denote by
$\langle\cdot,\cdot\rangle$ the canonical pairing between $M$ and
$N$. If we consider in $N_{\RR}$ the cone
 generated by the
elements of the standard basis $e_1,\ldots,e_n$, then the
corresponding toric variety is the affine space $\AA^n$ and the
subscheme $V(I)$ is invariant under the torus action (we refer to
\cite{Fu} for the basic notions on toric varieties). Hence the
normalized blowing-up of $\AA^n$ along $I$ is again a toric variety,
and therefore it corresponds to a fan subdividing the above cone.
This is the normal fan to the polyhedron $P_I$, that we will denote
by $\Delta_I$. It is defined as follows: to each face $Q$ of $P_I$
one associates the cone
$$\sigma_{Q}:=\{v\in N_{\RR}\vert\,\langle u,v\rangle\leq
\langle u',v\rangle\,{\rm for}\,{\rm every}
\,u\in Q\,{\rm and}\,u'\in P_I\}.$$ The
fan $\Delta_I$ consists of the cones $\sigma_Q$, when $Q$ varies
over the faces of $P_I$. Note that $\dim(\sigma_Q)=n-\dim(Q)$, so
the rays of $\Delta_{I}$ correspond to the facets of $P_I$, and
the maximal cones of $\Delta_{I}$ correspond to the vertices of
$P_I$.

Let $p$ be a fixed prime. We proceed now to the computation of
$Z_I$. For every $a=(a_1,\ldots,a_n)\in\NN^n$ we consider the set
$C_a=\prod_{i=1}^n(p^{a_i}\ZZ_p\smallsetminus p^{a_i+1}\ZZ_p)$.
Since each $p^{a_i}\ZZ_p\smallsetminus p^{a_i+1}\ZZ_p$ is a disjoint
union of $(p-1)$ translates of $p^{a_i+1}\ZZ_p$, we see that
$$\mu(C_a)=(p-1)^n\left(\frac{1}{p}\right)^{n+\sum_{i}a_i}.$$
We denote by $e$ the vector $(1,\ldots,1)$, so $\langle e,a\rangle
=\sum_ia_i$.

The function ${\rm ord}_pI$ is constant on $C_a$ with value
$$\nu(a):=\min\{\langle u,a\rangle\mid x^u\in I\}=\min\{\langle
u,a\rangle\vert u\in P_I\}.$$ Since the sets $C_a$ are disjoint and
the complement of their union has measure zero, we deduce
\begin{equation}\label{formula1}
Z_I(s)=\sum_{a\in\NN^n}\left(1-\frac{1}{p}\right)^n\cdot
\left(\frac{1}{p}\right)^{\langle e,a\rangle+s\nu(a)}.
\end{equation}
 Note that $\nu$ is a linear function on each of
the cones in $\Delta_{I}$. Indeed, if $w$ is a vertex of $P_I$, then
$\nu(a)=\langle w,a\rangle$ whenever $a$ is in $\sigma_w$.

If $\sigma$ is a cone in $\Delta_I$, choose a vertex $w$ of
$P_I$ such that $\sigma$ is contained in $\sigma_w$ and put
$\ell_{\sigma}:=w$. By letting the $a$ in (\ref{formula1}) vary
inside the relative interior of each cone in $\Delta_{I}$,
and using the definition in the Introduction, we get
the following

\begin{proposition}\label{formula2}
With the above notation, we have
\begin{equation}
Z_I(s)=\left(1-\frac{1}{p}\right)^n\cdot
\sum_{\sigma\in\Delta_{I}}Z_{\sigma,\ell_{\sigma},e}(s).
\end{equation}
\end{proposition}

\section{Igusa zeta functions for cones}

Our goal now in this section is to discuss Igusa-type zeta functions
associated to cones and prove Theorem~\ref{key}. Let $N$ be a
lattice, $M$ its dual, and
 $\sigma$  a pointed, rational polyhedral cone in $N_{\RR}$.
We consider $\ell_1$ and $\ell_2$ in $\sigma^{\vee}\cap M$, where
$\sigma^{\vee}$ is the dual cone of $\sigma$, such that $\sigma\cap
\{v\mid\ell_2(v)=0\}=\{0\}$. We want to study the function
$Z_{\sigma,\ell_1,\ell_2}$ and its poles.

The definition of $Z_{\sigma,\ell_1,\ell_2}$ was motivated by the
formula in Proposition~\ref{formula2}, but sometimes it is more
natural to consider a version of this function in which we sum over
all the integer points in $\sigma$:
\begin{equation}
\overline{Z}_{\sigma,\ell_1,\ell_2}(s):=\sum_{v\in\sigma\cap N}
\left(\frac{1}{p}\right)^{\ell_1(v)s+\ell_2(v)}.
\end{equation}
Again, this is well-defined if ${\rm Re}(s)>0$ and  we  have
$\overline{Z}_{\sigma,\ell_1,\ell_2}=\sum_{\tau}Z_{\tau,\ell_1,\ell_2}$,
where the sum is over all faces $\tau$ of $\sigma$. Moreover, it
follows from this formula that $Z_{\sigma,\ell_1,\ell_2}$ can be
computed in terms of the functions
$\overline{Z}_{\tau,\ell_1,\ell_2}$, where $\tau$ varies over the
faces of $\sigma$.

We start with the following

\begin{lemma}\label{lem1}
Let $v_1,\ldots,v_r$ in $N$ be linearly independent over $\QQ$. If
$w$ is in $N$ and $\ell_1$, $\ell_2$ are elements in $M$, with
$\ell_1$ nonnegative and $\ell_2$ positive on all the $v_i$, then we
put
$$A(s):=\sum_{v\in
S}\left(\frac{1}{p}\right)^{\ell_1(v)s+\ell_2(v)},$$ where
$S=\{w+a_1v_1+\ldots+a_rv_r\mid a=(a_i)\in\NN^r\}$. The function $A$
is well-defined and holomorphic for ${\rm Re}(s)> 0$ and it is a
rational function in $\left(\frac{1}{p}\right)^s$, so it has a
meromorphic continuation to $\CC$. Moreover, if $s$ is a pole of
$A$, then there is $i$ such that ${\rm
Re}(s)=-\frac{\ell_2(v_i)}{\ell_1(v_i)}$.
\end{lemma}

\begin{proof}
If ${\rm Re}(s)>-\frac{\ell_2(v_i)}{\ell_1(v_i)}$ for all $i$ such
that $\ell_1(v_i)$ is nonzero, then we have
$$A(s)=\left(\frac{1}{p}\right)^{\ell_1(w)s+\ell_2(w)}\cdot
\prod_{i=1}^n\sum_{a_i\in\NN}
\left(\frac{1}{p}\right)^{a_i(\ell_1(v_i)s+\ell_2(v_i))}$$
$$=\left(\frac{1}{p}\right)^{\ell_1(w)s+\ell_2(w)}\cdot
\prod_{i=1}^n
\frac{1}{1-\left(\frac{1}{p}\right)^{\ell_1(v_i)s+\ell_2(v_i)}}.$$
The assertions in the lemma are direct consequences of this formula.
\end{proof}

We can give now the proof of our result on Igusa-type zeta functions
associated to cones.

\begin{proof}[Proof of Theorem~\ref{key}]
Arguing by induction on the dimension of $\sigma$, we may assume
that the theorem holds for all cones of smaller dimension than
$\dim(\sigma)$ (the case when $\dim(\sigma)$ is zero being trivial).
In this case, we see that proving the assertions in the theorem for
$Z_{\sigma,\ell_1,\ell_2}$ is equivalent with proving them for
$\overline{Z}_{\sigma,\ell_1,\ell_2}$.

We show first that it is enough to prove the theorem when $\sigma$ is
a simplicial cone. Indeed, it is well-known that one can always find
a fan $\Gamma$ refining the cone $\sigma$ such that every cone in
$\Gamma$ is simplicial and the one-dimensional cones in $\Gamma$ are
precisely the rays of $\sigma$. Since
$$\overline{Z}_{\sigma,\ell_1,\ell_2}=\sum_{\tau\in\Gamma}
Z_{\tau,\ell_1,\ell_2},$$ and since each ray of a cone in $\Gamma$
is a ray of $\sigma$, we see that it is enough to prove the theorem
for each (maximal) cone in $\Gamma$.

Therefore we may assume that $\sigma$ is simplicial and our goal is
to show that $\overline{Z}_{\sigma,\ell_1,\ell_2}$ satisfies the
assertions in the theorem. Let $v_1,\ldots,v_r$ be the primitive
generators of the rays of $\sigma$. Since $\sigma$ is simplicial,
the $v_i$ are linearly independent over $\QQ$. The semigroup
$\sigma\cap N$ is finitely generated, so we may choose generators
$w_1,\ldots,w_s$. The $v_i$ span $\sigma$ over $\QQ$, hence we can
find a positive integer $m$ such that every $mw_j$ is in the
semigroup generated by the $v_i$. It follows that after replacing
$\{w_1,\ldots,w_s\}$ by $\{q_1w_1+\ldots+q_sw_s\mid 0\leq q_j\leq
m-1\}$, we may assume that
\begin{equation}\label{cover}
\sigma\cap N=\bigcup_{j=1}^s(w_j+S),
\end{equation}
 where $S$ is the
semigroup generated by the $v_i$.

If $I\subseteq\{1,\ldots,s\}$, let us put $$A_I(s):=\sum_v
\left(\frac{1}{p}\right)^{\ell_1(v)s+\ell_2(v)},$$ where the sum is
over $v$ in $\cap_{j\in I}(w_j+S)$. We claim that $\cap _{j\in I}
(w_j+S)$ is either empty or it is equal to $w+S$ for a suitable $w$
in $N$. Indeed, by an obvious induction on $|I|$ it is enough to
show this when $I$ has two elements,
say $j$ and $k$. The intersection of $w_{j}+S$ and
$w_{k}+S$ is nonempty if and only if $w_{j}-w_k$ lies in the lattice
generated by the $v_i$. If this is the case, let us write
$w_j-w_k=\sum_{i=1}^ra_iv_i$ for suitable integers $a_1,\ldots,a_r$.
If we put $w=w_j+\sum_{i=1}^r\max\{0,-a_i\}v_i$,
 then it is easy to check that $(w_j+S)\cap (w_k+S)=w+S$,
which proves our claim.

It follows from our claim and Lemma~\ref{lem1} that each $A_I$ is a
rational function of $\left(\frac{1}{p}\right)^s$. Moreover, if $s$
is a pole of $A_I$, then there is $i$ such that ${\rm
Re}(s)=-\frac{\ell_2(v_i)}{\ell_1(v_i)}$. On the other hand, it
follows from (\ref{cover}) that
\begin{equation}
\overline{Z}_{\sigma,\ell_1,\ell_2}=\sum_I(-1)^{|I|-1}A_I(s),
\end{equation}
where the sum is over all nonempty subsets $I$ of $\{1,\ldots,s\}$.
Therefore $\overline{Z}_{\sigma,\ell_1,\ell_2}$ satisfies the
assertions of the theorem, which completes the proof.
\end{proof}

Putting together Theorem~\ref{key} and the description of the Igusa
zeta function of a monomial ideal from the previous section we can
relate the poles of this zeta function with the  torus-invariant divisors
in the blowing-up along the ideal.

\begin{proof}[Proof of Theorem~\ref{thm2}]
It follows from Proposition~\ref{formula2} and Theorem~\ref{key}
that if $s$ is a pole of $Z_I$, then there is a primitive generator
$v$ of a ray of the normal fan $\Delta_{I}$ to the Newton
polyhedron $P_I$ such that
$${\rm Re}(s)=-\frac{\langle e,v\rangle}{\langle w,v\rangle}.$$
Here $w$ is a vertex of $P_I$ such that $v$ is contained in the
maximal cone $\sigma_w$ of $\Delta_I$ corresponding to $w$.

On the other hand, recall that the torus-invariant divisors on the
toric variety defined by $\Delta_{I}$ correspond precisely
to the rays of $\Delta_{I}$. Moreover, if $E$ is the divisor
corresponding to the ray through $v$, then $k_E=\langle e,v\rangle
-1$. Since we also have
$$a_E(I)=\min\{\langle u,v\rangle\mid x^u\in I\}=\langle
w,v\rangle,$$ as $v$ lies in $\sigma_w$, we deduce the assertion in
the theorem.
\end{proof}

\section{Poles and roots of the Bernstein-Sato polynomial}

We show now that the real part of any pole of $Z_I$ is a root of the
Bernstein-Sato polynomial $b_I$ associated to $I$. In fact, we prove
the following stronger statement that together with
Theorem~\ref{thm2} implies Theorem~\ref{thm1}.

\begin{proposition}\label{Bernstein}
If $I$ is a nonzero proper monomial ideal and if $E$ is a prime
divisor in the normalized blowing-up of the affine space along $I$
such that $a_E(I)$ is nonzero, then $-\frac{k_E+1}{a_E(I)}$ is a
root of the Bernstein-Sato polynomial $b_I$.
\end{proposition}

\begin{proof}
The divisor $E$ corresponds to a ray in the normal fan $\Delta_{I}$
to $P_I$. Let $v$ be a primitive generator of this ray. If $w$ is a
vertex of $P_I$ such that the corresponding maximal cone $\sigma_w$
of $\Delta_{I}$ contains $v$, then we have seen in the proof of
Theorem~\ref{thm2} that $k_E+1=\langle e,v\rangle$ and
$a_E(I)=\langle w,v\rangle$.
 Note that
since  $\langle w,v\rangle\neq 0$, the facet $Q$ of $P_I$
corresponding to $v$ is not contained in a coordinate hyperplane:
if, for example, $Q$ is contained in the hyperplane $(x_i=0)$, then
$v=e_i$ and since $w$ lies in $Q$ we get $\langle w,v\rangle=0$, a
contradiction.

In order to show that $(k_E+1)/a_E(I)$ is a root of the
Bernstein-Sato polynomial $b_I$ associated to $I$, we use the
description of the roots of $b_I$ from \cite{BMS2} (in fact, the
ones that we need for the theorem are ``the most straightforward''
of the roots of $b_I$). Since $Q$ is a facet of $P_I$ that is not
contained in a coordinate hyperplane, there is a unique linear
function $L_Q$ on $M_{\RR}$ such that $Q=P_I\cap L_Q^{-1}(1)$. With
this notation, it is shown in \cite{BMS2} (see Remark~4.6) that
$-L_Q(e)$ is a root of $b_I$.

On the other hand, since the ray through $v$ corresponds to the
facet $Q$ and since $w$ is in $Q$, we have
$$Q=\{u\in P_I\mid\langle u,v\rangle =\langle  w,v\rangle\}.$$
Therefore $L_Q$ is given by $\frac{1}{\langle w,v\rangle}\cdot v$
and since $-L_Q(e)$ is a root of $b_I$, we see that $(k_E+1)/a_E(I)$
 is, indeed,
a root of $b_I$.
\end{proof}

\begin{remark}
We do not know whether the analogue of Proposition~\ref{Bernstein}
holds for a non-necessarily monomial ideal $I$. Note that if $I=(f)$
is principal, then the assertion is trivial: the divisor $E$ is one
of the irreducible components of the divisor $H$ defined by $f$,
$k_E=0$ and $a_E(I)$ is the multiplicity of $E$ in $H$. The fact
that $-\frac{1}{a_E(f)}$ is a root of $b_f$ follows then by
restricting to an open subset where $E$ is smooth and $H=a_E(f)\cdot
E$.
\end{remark}

\begin{remark}\label{order}
The arguments in the previous two sections can be used to analyze
also the orders of the possible poles of the Igusa zeta function
$Z_I$. Indeed, it follows from Proposition~\ref{formula2} and from
the proof of Theorem~\ref{key} that if $s$ is a pole of order $r$ of
$Z_I$, then $r\leq n$ and there are $r$ invariant divisors
$E_1,\ldots,E_r$ on the normalized blowing-up along $I$ such that
$E_1\cap\ldots\cap E_r\neq\emptyset$ and ${\rm
Re}(s)=-(k_{E_i}+1)/a_{E_i}(I)$ for every $i$. We would like to
deduce that in this case ${\rm Re}(s)$ is a root of order $r$ of
$b_I$, but unfortunately, we do not understand well enough the
multiplicities of the roots of $b_I$.
\end{remark}

\begin{remark}\label{general_formula}
While Proposition~\ref{formula2} gives in principle a formula for
the Igusa zeta function of a monomial ideal, and Theorem~\ref{thm2}
gives an estimate on the denominator of this function (written as a
rational function of $1/p^s$), getting a general explicit formula
for the denominator is rather difficult. A \emph{Maple} code for
computing $p$-adic and motivic zeta functions of monomial ideals via
resolution of singularities is available, upon request, from the
first author.
\end{remark}

\begin{remark}\label{motivic}
Using motivic integration, Denef and Loeser defined in \cite{DL} a
motivic analogue of the Igusa zeta function. For concreteness, we
preferred to work with $p$-adic integrals. However,
 as the reader familiar
with this topic will certainly notice, all the above results have
analogues in the motivic setting, ``replacing $p$ by ${\mathbb
L}$''. For example, if $\sigma$, $\ell_1$ and $\ell_2$ are as in
Theorem~\ref{key}, then the series
\begin{equation}
\sum_{v\in\sigma\cap N}{\mathbb L}^{-(\ell_1(v)s+\ell_2(v))}
\end{equation}
can be written as a sum
of fractions with numerator in $K[{\mathbb L}^{-s}]$ and denominator
of the form
\begin{equation}
\prod_{i=1}^r\left(1-{\mathbb
L}^{-(\ell_1(v_i)s+\ell_2(v_i))}\right),
\end{equation}
 where $r\leq\dim(\sigma)$
and $v_1,\ldots,v_r$ are primitive generators of the rays of
$\sigma$. Here $K$ is the ring obtained from the Grothendieck ring
of varieties over a base field $k$ by inverting ${\mathbb
L}=[{\mathbb A}^1_k]$.
 Similarly, if $I$ is a monomial ideal, then the motivic zeta
 function of $I$
 \begin{equation}
\int_{({\mathbb A}^n)_{\infty}}{\mathbb L}^{-s\cdot{\rm ord}_tI}
 \end{equation}
can be written as a sum of fractions with numerator in $K[{\mathbb
L}^{-s}]$ and denominator of the form
$$\prod_{i=1}^r\left(1-{\mathbb L}^{-(a_{E_i}(I)s+k_{E_i}+1)}\right),$$
where $r\leq n$ and \textbf{}$E_1,\ldots,E_r$ are divisors on the
normalized blowing-up of ${\mathbb A}^n$ along $I$ such that
$E_1\cap\ldots\cap E_r$ is nonempty.
\end{remark}

\smallskip

\noindent{\bf Acknowledgements}. The second author would like to
thank Wim Veys for a very inspiring discussion on Igusa zeta
functions.

 \providecommand{\bysame}{\leavevmode \hbox
\o3em {\hrulefill}\thinspace}

\vfill\break

\end{document}